 \documentclass{article} 
 
\pdfoutput=1  

\usepackage{enumerate,wrapfig,placeins}

\usepackage{amsmath}
 
 \newtheorem{theorem}{Theorem}[section]
\newtheorem{lemma}[theorem]{Lemma}
\newtheorem{proposition}[theorem]{Proposition}

\newlength{\noteWidth}
\long\def\notes#1{\ifinner
             {\tiny #1}
             \else
              \marginpar{\parbox[t]{\noteWidth}{\raggedright\tiny #1}}
               \fi}

 \usepackage{textgreek,upgreek,bm}

 \usepackage{titlesec}

\usepackage[usenames, dvipsnames]{color}

\usepackage[pdfa,hidelinks]{hyperref}
  
\usepackage{optidef}
\usepackage{graphics} 
\usepackage{subfigure}

\graphicspath{{figures/}}

\usepackage{amsmath} 
\usepackage{amssymb}  
\usepackage{balance}

\usepackage{verbatim}

\usepackage{geometry} 

\def\urls#1{{\footnotesize\url{#1}}}

\def\mindex#1{\index{#1}}

\DeclareFontFamily{U}{mathx}{\hyphenchar\font45}
\DeclareFontShape{U}{mathx}{m}{n}{<-> mathx10}{}
\DeclareSymbolFont{mathx}{U}{mathx}{m}{n}
\DeclareMathAccent{\widebar}{0}{mathx}{"73}

\makeatletter
\newcommand\gobblepars{%
    \@ifnextchar\par%
        {\expandafter\gobblepars\@gobble}%
{}}
\makeatother

\def\whamrm#1{\smallbreak\pagebreak[3]%
	\noindent\text{\rm#1}\ \ \gobblepars}
	
\def\whamit#1{\smallbreak\pagebreak[3]%
	\noindent\textit{#1}\ \ \gobblepars}

\def\wham#1{\smallbreak\pagebreak[3]%
	\noindent\textup{\textbf{#1}}\ \ \gobblepars}

\def\Obj{\Upgamma}  

\def\Tdiff{\mathcal{D}}

\def\fee{\upphi}
\def\feex{\widetilde{\fee}}

\def\uH{\underline{H}}
\def\uQ{\underline{Q}}

\def\disc{\gamma}

\def\stepf{\beta}

\newcommand{\bbblot}{\raise1pt\hbox{\vrule height .4ex width .4ex depth .05ex}}

\long\def\defbox#1{\framebox[.9\hsize][c]{\parbox{.85\hsize}{%
\parindent=0pt
\baselineskip=12pt plus .1pt      
\parskip=6pt plus 1.5pt minus 1pt 
 #1}}}

\long\def\beginbox#1\endbox{\subsection*{}%
\hbox{\hspace{.05\hsize}\defbox{\medskip#1\bigskip}}%
\subsection*{}}

\def\endbox{}

 \def\archival#1{} 

\def\FRAC#1#2#3{\genfrac{}{}{}{#1}{#2}{#3}}

\def\ddtp{{\mathchoice{\FRAC{1}{d^{\hbox to 2pt{\rm\tiny +\hss}}}{dt}}%
{\FRAC{1}{d^{\hbox to 2pt{\rm\tiny +\hss}}}{dt}}%
{\FRAC{3}{d^{\hbox to 2pt{\rm\tiny +\hss}}}{dt}}%
{\FRAC{3}{d^{\hbox to 2pt{\rm\tiny +\hss}}}{dt}}}}

\def\ddyp{{\mathchoice{\FRAC{1}{d^{\hbox to 2pt{\rm\tiny +\hss}}}{dy}}%
{\FRAC{1}{d^{\hbox to 2pt{\rm\tiny +\hss}}}{dy}}%
{\FRAC{3}{d^{\hbox to 2pt{\rm\tiny +\hss}}}{dy}}%
{\FRAC{3}{d^{\hbox to 2pt{\rm\tiny +\hss}}}{dy}}}}

\def\argmin{\mathop{\rm arg{\,}min}}

\def\state{{\sf X}}

\def\ustate{{\sf U}} 
\def\ystate{{\sf Y}}

\def\ystate{{\sf Y}}

\def\bfmath#1{{\mathchoice{\mbox{\boldmath$#1$}}%
{\mbox{\boldmath$#1$}}%
{\mbox{\boldmath$\scriptstyle#1$}}%
{\mbox{\boldmath$\scriptscriptstyle#1$}}}}

\def\bfPhi{\bfmath{\Phi}}

\def\bfmX{\bfmath{X}}

\def\bfmY{\bfmath{Y}}

\def\bfmhhaY{\bfmath{\hhaY}} 
\def\bfmhhaY{\hbox to 0pt{$\widehat{\bfmY}$\hss}\widehat{\phantom{\raise 1.25pt\hbox{$\bfmY$}}}}

\def\hatheta{{\hat\theta}}

\def\haR{\widehat R}

\def\clC{{\cal C}}

\def\clS{{\cal S}}

\def\clX{{\cal X}}
\def\clY{{\cal Y}}

\def\clC{{\cal C}}

\def\eqdef{\mathbin{:=}}

\def\Prob{{\sf P}}

\def\Expect{{\sf E}}

\def\Cov{\hbox{\sf Cov}}

\def\ind{\bbbone}
  
 \def\epsy{\varepsilon}

\def\formtmp#1#2{{\vskip12pt\noindent\fboxsep=0pt\colorbox{#1}{\vbox{\vskip3pt\hbox to \textwidth{\hskip3pt\vbox{\raggedright\noindent\textbf{#2\vphantom{Qy}}}\hfill}\vspace*{3pt}}}\par\vskip2pt%
\noindent\kern0pt}}

\def\barf{{\widebar{f}}}

\def\barA{{\bar{A}}}

\def\barJ{{\bar{J}}}

{\end{list}}

\def\ass(#1:#2){(#1\ref{#1:#2})}

\def\ritem#1{
\item[{\sf \ass(\current_model:#1)}]
}

\newenvironment{recall-ass}[1]{%
\begin{description}
\def\current_model{#1}}{
\end{description}
}

 \newcommand{\blot}{\vrule height 1.1ex width .9ex depth -.1ex }
\def\qedb{\ifmmode\blot\else{\vspace{-.2cm}\unskip\nobreak\hfil
\penalty50\hskip1em\null\nobreak\hfil\blot
\parfillskip=0pt\finalhyphendemerits=0\endgraf}\fi}

  \makeatletter
\DeclareRobustCommand{\sqcdot}{\mathbin{\mathpalette\morphic@sqcdot\relax}}
\newcommand{\morphic@sqcdot}[2]{%
  \sbox\z@{$\m@th#1\centerdot$}%
  \ht\z@=.33333\ht\z@
  \vcenter{\box\z@}%
}
\makeatother

\newcommand{\qedsymbol}{\hbox{\tiny$\blacksquare$}}  %

\def\qed{\ifmmode\qedsymbol\else{\unskip\nobreak\hfil
\penalty50\hskip1em\null\nobreak\hfil\qedsymbol
\parfillskip=0pt\finalhyphendemerits=0\endgraf}\fi}

\newcounter{rmnum}

\newcounter{anum}
\newenvironment{alphanum}{\begin{list}{{\upshape (\alph{anum})}}{\usecounter{anum}
\setlength{\leftmargin}{18pt}
\setlength{\rightmargin}{8pt}
\setlength{\itemindent}{2pt}
}}{\end{list}}

\newcommand{\field}[1]{\mathbb{#1}}

\def\Re{\field{R}}

\def\Co{\field{C}}

\def\Prob{{\sf P}}

\def\Expect{{\sf E}}

\def\transpose{{\intercal}}

\def\argmin{\mathop{\rm arg\, min}}
\def\ind{\hbox{\large \bf 1}}

\def\epsy{\varepsilon}

\def\haJ{\widehat J}

\def\haY{\widehat{Y}}

\def\hhaY{\hbox to 0pt{$\haY$\hss}\widehat{\phantom{\raise 1.25pt\hbox{Y}}}}

\def\haY{\widehat Y}

\def\bfPhi{\bfmath{\Phi}}

\newlength{\dhatheight}

\def\EpsLength{\upxi}
\def\hafee{\widehat\fee}

\def\whamb{\wham{$\bullet$} }

\def\thetaPR{\theta^{\text{\tiny\sf  PR}}}
\def\thetaPR{\theta^{\text{\tiny\sf  PR}}}
\def\tilthetaPR{\tilde{\theta}^{\text{\tiny\sf  PR}}}
\def\barthetaPR{\bar{\theta}^{\text{\tiny\sf  PR}}}

\def\tTheta{{\text{\tiny$\Theta$}}}

\def\thetaPR{\theta^{\text{\tiny\sf  PR}}}
\def\tilthetaPR{\tilde{\theta}^{\text{\tiny\sf  PR}}}
\def\SigmaPR{\Sigma^{\text{\tiny\sf PR}}_\tTheta}

\usepackage[usenames,dvipsnames]{color}
\usepackage{color}
 \definecolor{programcode}{gray}{0.9}
 
 \definecolor{lightgray}{gray}{0.7}

\definecolor{MyDarkBlue}{cmyk}{0.5,0.1,0,0.9}

\def\bl#1{{\color{blue}#1}}

\usepackage{cleveref}

\Crefname{corollary}{Corollary}{Corollaries}
\Crefname{eqnarray}{eq.}{eqs.}
\Crefname{equation}{eq.}{eqs.}

\Crefname{figure}{Fig.}{Figs.}
\Crefname{tabular}{Tab.}{Tabs.}
\Crefname{table}{Tab.}{Tabs.}
\Crefname{lemma}{Lemma}{Lemmas}

\Crefname{theorem}{Thm.}{Thms.}
\Crefname{definition}{Definition}{Definitions}
\Crefname{section}{Section}{Sections}
\Crefname{proposition}{Prop.}{Propositions}
\Crefname{assumption}{Assumption}{Assumptions}
\Crefname{example}{Example}{Examples}

\def\whamit#1{\smallbreak\pagebreak[3]%
	\noindent\textit{#1}\ \ \gobblepars}

\def\wham#1{\smallbreak\pagebreak[3]%
	\noindent\textbf{#1}\ \ \gobblepars}

\def\bl#1{{\color{blue}#1}}

\def\logderfeex{\Uplambda} 
\def\SGDnabla{{\breve{\nabla}\!}}

\def\sstate{\textsf{S}}
\def\preDens{f_0}
\def\postDens{f_1}

\def\pFA{p_{\textsf{\tiny FA}}}

\def\Obs{Y}
\def\bfObs{\bm{Y}}

\def\postObs{X^1}
\def\bfpostObs{\bfmX^1}
\def\preObs{X^0}
\def\bfpreObs{\bfmX^0}

\def\condDist{\Uppi} 

\def\InfoState{\clX}    
\def\surf{\breve{f}}

\def\thexp{\upupsilon}

\def\tchange{\uptau_{\text{\scriptsize\sf a}}}

\def\tstop{\uptau_{\text{\scriptsize\sf s}}}

\def\expa{\varrho_{\text{\scriptsize\sf a}}}

 \newcommand{\overbar}[1]{\mkern 1.5mu\overline{\mkern-1.5mu#1\mkern-1.5mu}\mkern 1.5mu}

\def\thresh{%
	\mathchoice
	{\text{\small\rm H}}%
	{\text{\small\rm H}}%
	{\text{\scriptsize\rm H}}%
	{\text{\tiny\rm H}}}

\def\barthresh{\overbar{\thresh}}
\def\hathresh{\widehat{\thresh}}

\def\hathresh{\widehat{\thresh}}

\def\invTemp{\upxi}

\def\MDD{{\sf MDD}}
\def\MDE{{\sf MDE}}

\def\CauchyScale{\gamma}

 \def\bqed{{\color{blue} \qedb} \bigskip}

\title{Reinforcement Learning Design for Quickest Change Detection}

\author{Austin Cooper
	\and Sean Meyn
	\thanks{ASC and SPM are with the University of Florida, Gainesville, FL 32611 Financial support from ARO award W911NF2010055     and NSF award  CCF 2306023
		   is gratefully acknowledged.}%
}

\begin{document}
 
\maketitle

 \begin{abstract}
The field of quickest change detection (QCD) concerns design and analysis of algorithms to estimate in real time the   time at which an important event takes place,   
and identify properties of the post-change behavior.
  
It is shown in this paper that approaches based on reinforcement learning (RL) can be adapted based on any ``surrogate information state'' that is adapted to the observations.     Hence we are left to choose both the surrogate information state process and the algorithm.  For the former, it is argued that there are many choices available, based on a rich theory of asymptotic statistics for QCD.  
 Two approaches to RL design are considered:
\begin{alphanum}
\item 
Stochastic gradient descent based on an actor-critic formulation.  Theory is largely complete for this approach:  the algorithm is unbiased, and will converge to a local minimum.  However, it is shown that variance of stochastic gradients can be very large, necessitating the need for commensurately long run times. 

\item
Q-learning algorithms based on a version of the projected Bellman equation. 
It is shown that the algorithm is stable,  in the sense of bounded sample paths, and that a solution to the projected Bellman equation exists under mild conditions. 

\end{alphanum}
Numerical experiments illustrate these findings, and provide a roadmap for algorithm design in more general settings.

\end{abstract}

\tableofcontents

\clearpage
\section{Introduction} 
\label{s:intro}

The goal of the
research surveyed in this paper is to create algorithms for \textit{quickest change detection} (QCD),
 for applications in which statistics are only partially known, particularly after the change has occurred.   
 While the authors were initially motivated by applications in power systems, the setting here is entirely general.   
 Examples of events that we wish to detect include  human or robotic intruders, computer attack, faults in a power system,  and onset of heart attack for a patient~\cite{liatarvee21,liavee22}.  
  
The standard QCD model includes a sequence of observations $\bfObs\eqdef \{ \Obs_k :  k\ge 0\}$, assumed here to evolve as a real-valued  stochastic process. 
The statistics of these observations change at a time denoted $\tchange\ge 0$. 
The goal is to construct an estimate of the change time, denoted $\tstop$, that is adapted to the observations.  That is, on denoting $\Obs_0^k = (\Obs_0; \cdots; \Obs_k)$,  for each $k$ we may write $\ind\{\tstop \le k\} =  s_k(\Obs_0^k)$ for some Borel-measurable mapping $s_k\colon \Re^{k+1} \to \{0,1\}$. 
The estimate must balance  two costs: 
1.  \textit{Delay}, which is expressed $(\tstop -\tchange)_+\eqdef \max(0,\tstop -\tchange)$,  and 2.\ \textit{false alarm}, meaning that $\tstop -\tchange<0$.

There are two general models that lead to practical solutions:  
Bayesian and minimax approaches.
Typical measures of performance for the former approach are based on \textit{mean detection delay} $\MDD$ and probability of false alarm $\pFA$: 
 \begin{equation}
\text{ $\MDD  =   \Expect[ (\tstop -\tchange)_+ ]  $  
		\ \ 
			and 
		\ \ 
$\pFA=  \Prob\{\tstop < \tchange \}  $.}
\label{e:MDDpFA}
\end{equation}

The focus of this paper is on the Bayesian approach,  based on a partially observed Markov Decision Process (POMDP).   See \Cref{s:POMDP} for canonical examples.

%

Successful approaches to algorithm design are typically based on the construction of a real-valued stochastic process $\{\InfoState_n\} $ that plays a role similar to the celebrated information state of POMDP theory, 
and a threshold policy is adopted: for a pre-assigned threshold $\thresh>0$, the stopping rule is    
\begin{equation}
\tstop = \min\{ n\ge 0 :   \InfoState_n\ge \thresh \}\, .
\label{e:threshold}
\end{equation}%
\begin{subequations}%
Two famous examples are  defined recursively:  with $\InfoState_0=0$,  
\begin{align}
&\text{1. Shiryaev–Roberts:}   &&
\InfoState_{n+1} =   \exp\bigl( L_{n+1} \bigr) [\InfoState_n  +  1] 
\label{e:SR}
	\\
&\text{2. CUSUM:}   &&
\InfoState_{n+1} =    \max\{ 0,  \InfoState_n +  L_{n+1}  \} 
\label{e:CUSUM}
\end{align}
in which $L_n= L(Y_n)$ is a log likelihood ratio for the conditional i.i.d.\ settings in which these models are typically posed (see \Cref{s:BQCD}). In particular, the CUSUM statistic evolves as a reflected random walk (RRW) with negative drift  for $0\le n < \tchange$.
\label{e:2surInfoState}
\end{subequations}

Analysis of the threshold policy \eqref{e:threshold} is typically posed in an asymptotic setting, considering a sequence of models with threshold $\thresh$ tending to infinity.    Approximate optimality results for either statistic may be found in   \cite{mou86,mostar09}.
See \cite{poltar10,xiezouxievee21}  for further  history.

\wham{Contributions}

This paper develops theory for QCD in a Bayesian setting,  and demonstrates how the solution structure lends itself to RL design.   One theme is the application of observation-driven statistics such as \eqref{e:CUSUM} to form a ``surrogate'' information state for policy synthesis.  

\whamb  
Performance of the CUSUM test is approximated in the asymptotic setting in which there is a strong penalty for false alarm. Well-known analysis for the ideal case is extended to settings for which $L_n$ in \eqref{e:CUSUM} is mismatched with the observations (i.e., $L_n$ is not a log likelihood ratio). 
This general theory motivates the control architectures proposed for RL design.    In particular: 

\wham{$\triangleright$}
 An actor-critic approach is introduced and shown to be consistent under mild conditions (see \Cref{t:SGDconverges}).
 
\wham{$\triangleright$}
 A Q-learning algorithm is introduced and shown to be stable provided the input used for training is \textit{sufficiently optimistic} \cite{mey23}.

 \whamb    The theory is illustrated with many experiments, comparing resulting policies with common heuristics as well as the true optimal.  Among the findings are
 \wham{$\triangleright$}
 Stability for the scalar gain algorithm requires extremely high level of optimism, resulting in poor numerical performance.     A version of Zap Q-learning is far more reliable.     
   \wham{$\triangleright$}
  The resulting policies performed well using a basis obtained via binning, and a linear function class inspired by results obtained via binning.

\wham{Literature}

See \cite{liatarvee21,liavee22} for excellent recent surveys on QCD theory.   Much of this theory is cast in a minimax rather than Bayesian setting.  Numerical techniques  to solve the QCD problem in the Bayesian setting may be found in \cite{zhoenl13}.

The analysis in \Cref{s:asystats} is cast 
as 
the conditionally i.i.d.\  model of \cite{shi77},  though we relax the classical assumption that the change time has a geometric distribution.

The vast majority of theory requires statistical independence of $\bfpreObs$, $\bfpostObs$ and $\tchange$.   This is the case  in Shiryaev's conditional i.i.d.\ model,  for which the   stochastic processes $\bfpreObs$, $\bfpostObs$ are also assumed i.i.d.;

Extension to a conditionally Markov model or hidden Markov model is possible by adapting techniques from the recent work \cite{zhasunherzou23,zhasunherzou23b};  while cast in an adversarial setting, many approximations remain valuable in the Bayesian setting adopted here.

Stability theory of  Q-learning  for optimal stopping was resolved in \cite{tsivan99};   conditions for consistency are similar to those for the simpler TD-learning algorithm.   However, the specific algorithm considered required that the cost function be fully observed.  This is why Q-learning is re-considered in the present article.  

In this prior work it is recognized that the state is inherently partially observed.    In   \cite{tsivan99} along with many papers in the RL literature,   
a truncated history of observations is adopted as a surrogate information state,
$\InfoState_k = (Y_{k-B+1};\cdots;Y_k)$, with $B\ge 1$. An innovations process obtained from the Kalman filter is used to define $\{\InfoState_k \}$ in applications to power systems   \cite{kuroguliwan19}.
General theory surrounding the approximation of the information state may be found in    \cite{subsinsermah22} (along with substantial history).

There is a long history of application of techniques from 
RL to approximate the solution to the optimal stopping problem.  
The first stability analysis of Q-learning with linear function approximation appeared in \cite{tsivan99},  which inspired significant research such as \cite{liszesch09,chedevbusmey20a}.    
The algorithms conceived in this prior work are not applicable in the applications considered in this paper because 
the cost (or rewards) is assumed to be fully observed.   
The RL algorithms introduced in this paper are more complex, and have a weaker supporting stability theory, precisely because this assumption is violated.

\wham{Organization}  

\Cref{s:BQCD} provides background on the standard Bayesian QCD problem as well as an alternate cost criterion for which approximations are formulated.
\Cref{s:RLQCD} includes formulations of two RL approaches to optimal stopping. The paper then turns to design and experimental findings of Q-learning applied to our Bayesian QCD problem in \Cref{s:Qqcd}. \Cref{s:conc} provides concluding thoughts and directions for future work. 
This preprint is an extension of \cite{coomey24d}, to appear in the 2024 IEEE Conference on Decision and Control.

\section{Bayesian QCD}
\label{s:BQCD}

This section contains background on approaches to modeling and algorithm design for QCD.   We begin with a canonical Bayesian model, cast as a POMDP.

\subsection{POMDP model}
\label{s:POMDP}

In this model both the change time and the observations are deterministic functions of a   time-homogeneous Markov chain $\bfPhi$, evolving on a state space $\state$.   It is assumed that $\Obs_k = h(\Phi_k)$,  $k\ge 0$,  for a function $h\colon\state\to\ystate$ (measurable in an appropriate sense).   
Assume moreover that there is a decomposition
$\state = \state_0 \cup \state_1$,   for which $\state_1$ is \textit{absorbing}:   $\Phi_k\in \state_1$ for all $k\ge 0$ if $\Phi_0\in \state_1$.     
The change time is defined by $\tchange = \min\{k\ge 0 :  \Phi_k \in \state_1 \}$.

We arrive at a POMDP with input $U_k \in  \ustate = \{0,1\}$, and  $\tstop$ defined as the first value of $k$ such that $U_k=1$.
The control problems of interest are optimal stopping problems:     
 For any cost functions $c_\circ ,c_\bullet \colon\state    \to \Re$, we wish to minimize  over all inputs adapted to the observations,
\begin{equation}
 J(\Phi_0, U_0^\infty ) =   \Expect\Big[ \sum_{k=0}^{\tstop -1}  c_\circ(\Phi_k  )     + c _\bullet(\Phi_{\tstop}  )  \Big]
\label{e:ObjPOMDPQCD}
\end{equation}
Consistent with the standard QCD framework is $c_\circ(z) =  \ind\{ z \in \state_1\}  $  and $c_\bullet(z) = \kappa  \ind\{ z\in \state_0 \}$  with $\kappa>0$,  so that $ J(\Phi_0, U_0^\infty )  = \MDD + \kappa \pFA$  (recall \eqref{e:MDDpFA}).

The structure of an optimal solution can be expressed as state feedback with suitable choice of state process.  We use the term \textit{information state},  denoted $\{\InfoState_k : k\ge 0\}$.  This is defined as a sufficient statistic for optimal control, in the sense that    an optimal solution is expressed as ``information state feedback'',  $U_k^* = \fee^*(\InfoState_k)$.    The canonical example is 
 $\{\InfoState_k \} =  \{ \condDist_k \}$,  the sequence of conditional distributions (often called the belief state) \cite{elllakmoo95,kri15}.   This structure leads to a practical solution when $\state$ is finite, with $K$ elements,  so that $\condDist_k$ evolves on the $K$-dimensional simplex $\clS^K$,   and $ \fee^* \colon \clS^K \to  \ustate$ is measurable.

\wham{Shiryaev's model}

The POMDP model is a generalization of Shiryaev's conditional i.i.d.\ model, in which observations are expressed
\begin{equation}
\Obs_k = \preObs_k \ind_{k< \tchange } +  \postObs_k \ind_{k \ge \tchange }  \,,\qquad k\ge 0\,, 
\label{e:QCDmodel}
\end{equation} 
with   $\bfpreObs$ and $\bfpostObs$ i.i.d.\  and mutually independent stochastic processes;
the change time  $\tchange$ is independent of  $\bfpreObs, \bfpostObs$, and has a geometric distribution.   Under these strong assumptions,  the real-valued process $\{ p_k =  \Prob\{  \tchange \le  k \mid \Obs_0^k \}  : k\ge 0\}$   serves as an information state,    and an optimal test is of the form $U_k^* = \ind\{  p_k \ge \thresh\} $ for some   threshold $\thresh>0$ (see \cite{shi77} and the tutorial \cite{veeban14}).

The observation model \eqref{e:QCDmodel} is valuable in analysis of common heuristics.   
Suppose that the marginal distributions of $ \{\preObs_k  ,  \postObs_k\}  $  have densities on $\Re$, denoted  $ \preDens $, $\postDens$, and denote $L(y)=\log(\postDens(y)/\preDens(y)) $---the log likelihood ratio (LLR).  
Crucial for analysis of either of the algorithms \eqref{e:2surInfoState} is that $L$ has  positive mean under $\postDens$ and negative mean under $\preDens$.    

It is known that either of the algorithms \eqref{e:2surInfoState}  is approximately optimal for Shiryaev's model,  
for large $\kappa$  and large mean change time  $ \Expect[\tchange]$  
\cite{xiezouxievee21}.

\wham{Alternative to the standard cost criterion}

The standard cost criterion is $\MDD + \kappa \pFA$   is sensible in Shiryaev's model in which the change time is independent of 
$\{ \preObs_k ,  \postObs_k :  k\ge 0\}$.  In the general POMDP model there may be evidence that a change is imminent;   in such cases, a (common sense) good decision rule might make an early declaration of change. 
These decision rules might be far from optimal under the usual cost criterion since it is insensitive to the value of
\textit{eagerness}, defined as  $(\tstop -\tchange)_- \eqdef \max(0, - [\tstop -\tchange] )$.    In this paper we consider  
the \textit{mean detection eagerness} $\MDE=\Expect[ (\tstop -\tchange)_- ]$ in the cost criterion $\MDD+\kappa\MDE$, leading to what we believe is a more reasonable objective,
 \begin{equation}
 J(\Phi_0, U_0^\infty ) =   \Expect\big[ (\tstop -\tchange)_+  +  \kappa (\tstop -\tchange)_-  \big]
\label{e:EagerObj}
\end{equation}
This may be placed in the POMDP standard form \eqref{e:ObjPOMDPQCD}, with
\begin{equation}
c_\circ(z) =  \ind\{ z \in \state_1\} \,, \qquad
c_\bullet(z) = \kappa \Expect[  \tchange  \mid \Phi_0 = z] 
\label{e:eager}
\end{equation}

\subsection{Asymptotic statistics}
\label{s:asystats}

The remainder of this section concerns the CUSUM test. Analysis is restricted to the following conditional i.i.d.\ setting:   The stochastic processes 
$\{  \preObs_k \}$,   $\{  \postObs_k \}$  are each i.i.d.,   and independent of the change time $  \tchange$. 
It is assumed that the marginal densities $f_0$ and $f_1$ exist,  and that the LLR  $L = \log(f_1/f_0)$ exists  and  is integrable with respect to either $f_1$ or $f_0$   

Our interest is approximating the performance of the CUSUM test, and also approximating the optimal threshold for a given value of $\kappa$. 
The analysis allows for two significant relaxations:
\wham{1.}  
We consider  $L_n = F( \Obs_n) $ for a Borel measurable function $ F\colon \ystate \to\Re $, not necessarily the LLR.     Letting $m_i = \int F(y) \,  f_i(y) dy$ for $i=0,1$, it is assumed that $m_0<0$ and $m_1>0$.    Hence  the RRW \eqref{e:CUSUM} is a positive recurrent Markov chain if $\tchange =\infty$.

Two log moment generating functions  are denoted   $\Lambda_i(\thexp) = \log \Expect[\exp(\thexp F(X^i_k)) $  for $\thexp\in\Re$ and $i = 0,1$.

\wham{2.}  
The strong distributional assumption on the change time is replaced by a regularity condition:    
\wham{Regular geometric tail:} for some $\expa < \infty$, 
\begin{equation}
	\lim_{n\to\infty} \frac{1}{n} \log  \Prob\{ \tchange \ge   n  \} = -\expa 
	\label{e:hazardAss}
\end{equation}
 
The regularity assumption obviously holds in Shiryaev's model, in which $\tchange$ has a  geometric distribution.    We obtain $\expa>0$ in the POMDP model under mild assumptions.  
Full proofs for results in this section may be found in \cite{coomey24b,coomey24c}, along with extensions beyond the conditionally i.i.d.\ model.

\begin{lemma}
\label[lemma]{t:POMDPhazard}
Consider the POMDP model  with  $\state$   finite, 
$\Prob\{ \tchange  = \infty \} = 0$,  yet 
 $\Prob\{ \tchange  > N \} >0$ for each  $N>0$ and  $\Phi_0 \in\state_1$.
       Then  \eqref{e:hazardAss} holds for some   $\expa>0$.
\end{lemma}

The cost of delay is easily approximated for this model:   After a change has occurred, the most likely path is linear with slope $m_1 >0$.   For a threshold $\thresh\gg 1$,  the delay  $ (\tstop -\tchange)_+$  is overwhelmingly likely to be close to $\thresh/m_1$.   

Approximation of the mean of $ (\tstop -\tchange)_-$ is based on  
well-established large deviations theory for RRWs. The main results of this theory require that the  log moment generating functions $\Lambda_i (\thexp) \eqdef \log \int \exp(\thexp x) f_i(x)\, dx$,  $i=0,1$,  be finite over a suitable range of $\thexp\in\Re$.

\begin{subequations}
	
Denote
for any threshold test,  
\begin{align}  
	\thresh^*(\kappa) &= \argmin_{\thresh\ge 0} \{ \kappa \MDE(\thresh) +  \MDD(\thresh)  \} 
	\label{e:optThreshold}
	\\
	J^*(\kappa) &= \min_{\thresh\ge 0} \{ \kappa \MDE(\thresh) +  \MDD(\thresh)  \}   
	\label{e:optCost}
\end{align}%
\label{e:optThresholdCost}%
\end{subequations}%
We write CUSUM* to denote the CUSUM algorithm using the optimal threshold $\thresh^*(\kappa)  $. In \Cref{t:barcApprox} we justify these approximations,
\begin{equation}
	\barthresh_\infty^*( \kappa)  = \frac{1}{  \thexp_+ }          \log(\kappa)    
	\,,	\quad 
	\barJ_\infty^*( \kappa) =
	\frac{1}{  m_1 }   \frac{1}{  \thexp_+ }        \log(\kappa) 
	\label{e:ApproxJH}
\end{equation}
where the constant $\thexp_+$ is defined in the proposition.

\begin{subequations}
	
	\begin{proposition}
		\label[proposition]{t:barcApprox}
		Suppose the following conditions hold:   
\whamrm{1)}  the limit \eqref{e:hazardAss}  holds with $\expa >0$,    
\whamrm{2)}  $\Lambda_0$ has  two distinct roots $\{0,\thexp_0 \}$,   a unique solution 
		$\thexp_+ > \thexp_0$  to     $\Lambda_0(\thexp_+) = \expa$, and $\Lambda_0$ is   finite-valued in a neighborhood of $[0,\thexp_+]$,   
\whamrm{3)}  $\Lambda_1$  is finite-valued in a neighborhood of the origin. 

		Then, 
		\begin{equation}
			\barJ (\thresh, \kappa)  = \thresh  [  1 /m_1  +  o(1)  			
			] +  \kappa \exp(-\thresh[  \thexp_+  +  o(1) 			 			 ] )
			\label{e:Japprox1}
		\end{equation}
		in which  $ o(1) \to 0$  as $\thresh\to\infty$ for $i=1,2$.  
		Consequently,   for the CUSUM* test,
		\begin{align}
\thresh^*(\kappa) & =  \barthresh_\infty^*( \kappa )	  +   o(\log(\kappa))
\label{e:Japprox2} 
\\
J^*(\kappa) 
& =  \barJ_\infty^*( \kappa) 
+   o(\log(\kappa)) 
%
\label{e:Japprox3}
		\end{align}
		\qed
	\end{proposition}
	
\end{subequations}

\begin{figure}[h]
	\centering		
	\includegraphics[width=0.65\hsize]{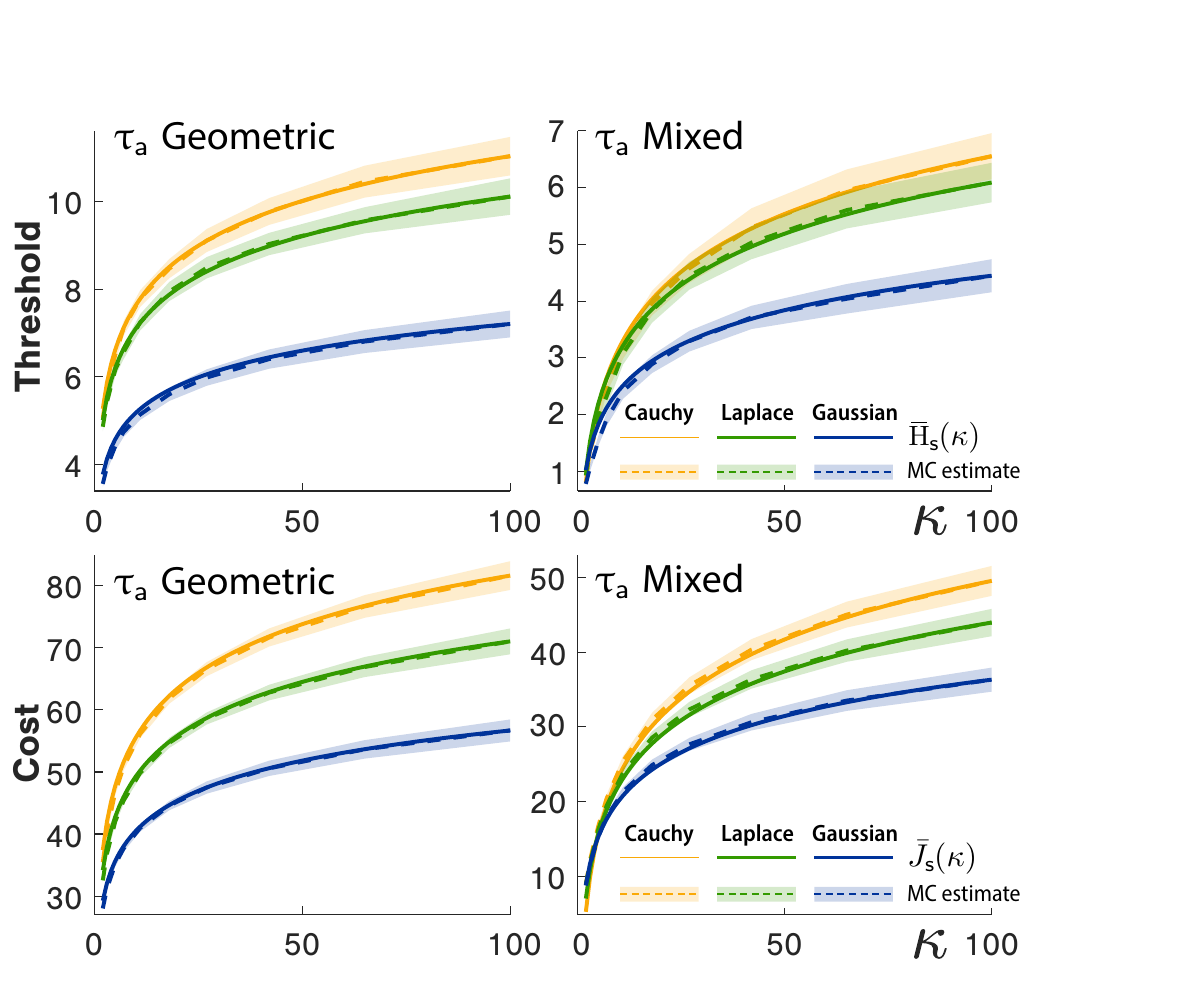} 
	\caption{Approximations of the optimal threshold and cost.}
	\label[figure]{f:CostAndThresholdApproximations}
\end{figure}

Approximations for three choices of $F$ are shown in \Cref{f:CostAndThresholdApproximations}. Details can be found in \Cref{s:Qqcd}. 

\section{Reinforcement Learning and QCD}
\label{s:RLQCD}

In the few examples we have considered we have found that the approximations in \Cref{t:barcApprox} are highly accurate.   
In non-ideal settings the proposition is valuable in the construction of RL algorithms.   We provide algorithms, and full justification in some cases.   Algorithm design and analysis is set in the POMDP (Bayesian) setting,   with  cost criterion \eqref{e:EagerObj}.

Assumed given is a surrogate belief state:   a stochastic process $ \{\InfoState_k : k\ge 0\}$, evolving on a closed subset of Euclidean space $\sstate$, and which is adapted to the observations $\clY_k =\sigma\{ Y_0,\dots, Y_k\}$.  
 We do not require that $\InfoState_k$ is in any sense an approximation of an information state.  
  In particular,  the numerical experiments largely focus on 
the CUSUM statistic   \eqref{e:CUSUM}.

We first consider a version of the actor-critic method,  followed by approaches to Q-learning.   

\wham{Actor-critic method}

Assumed given is a collection of randomized stationary policies $\{ \feex^\theta : \theta \in \Re^d \}$.
For each $\theta$ the statistics of the decision rule are defined  via
\[
\Prob\{U_k=u \mid \Obs_0^k;   \InfoState_k = x\}   =  \feex^\theta ( u \mid x)   \,,    \ \  u\in\ustate\,, \ x\in\sstate   .
\]

We fix  an initial distribution $\nu$ for the   Markov chain $\Psi$, and  denote $\mu_\theta (z , u) = \nu(z)   \feex^\theta ( u \mid x) $ for $z = (x; \varsigma) \in \state\times\sstate$  and $u\in \ustate$.    Our goal is to minimize  
\begin{equation}
\Obj(\theta) =  \Expect^\theta_{\mu_\theta}\Big[ \sum_{k=0}^{\tstop }  c(\Phi_k, U_k  )     \Big]
\label{e:Obj}
\end{equation}
The subscript   indicates that $ (\Phi_0,   \InfoState_0 ,   U_0 ) \sim \mu_\theta$,  and the superscript ``$\theta$'' indicates that  the policy $\feex^\theta$ determines the input.

We consider   stochastic gradient descent (SGD),
\begin{equation}
\theta_{n+1} = \theta_n - \alpha_{n+1} G_n \SGDnabla_\Obj(n)\,,
\label{e:QCDAC}
\end{equation}
 in which the
stepsize $  \alpha_{n+1} $ and
 the matrix gain $G_n$ are   design choices.

In the actor-critic algorithm the stochastic   gradient $ \SGDnabla_\Obj(n)$ is represented    in terms of the \textit{score function},
\begin{equation}
\logderfeex^\theta(x,u) =  \nabla_\theta \log[\feex^\theta(u\mid x) ]  
\label{e:ScoreFunction}
\end{equation} 
  defined to be zero for any values for which $\feex^\theta(u\mid x) =0$.

The following result follows from a long history surveyed in the Notes section of \cite[Ch. 10]{CSRL}:

\begin{subequations}
	
	\begin{proposition}
		\label[proposition]{t:ScoreRep}
		Suppose that $ \Obj$  and $\nabla\Obj$  are continuous.   Then  $\nabla\Obj\, (\theta) =   \Expect^\theta_{\mu_\theta}[ \SGDnabla_\Obj  ]$,  for either of the two options:
		\begin{align}
			\SGDnabla_\Obj^\theta  &=    \sum_{k=0}^{\tstop }  c(\Phi_k, U_k  ) S_k^\theta 
			\label{e:nablaObjElig}
			\\
			\textit{or}   \qquad
			\SGDnabla_\Obj^\theta & =    \sum_{k=0}^{\tstop }  
			Q_\theta(\Psi_k, U_k  )  \logderfeex^\theta_k
			\label{e:nablaObjQ}
		\end{align} 
		in which $ \logderfeex^\theta_k \eqdef \logderfeex^\theta(  \Phi_k , U_k)  $, $S_k^\theta = \logderfeex^\theta_0+\cdots+ \logderfeex^\theta_k$, and $Q_\theta$ is the fixed-policy score function 
		\begin{equation}
			Q_\theta(z,u   )  =   \Expect^\theta \Big[ \sum_{k=0}^{\tstop }  c(\Phi_k, U_k  ) \mid \Psi_0=z, U_0=u   \Big]
			\label{e:QAC}
		\end{equation}
	\end{proposition}
\end{subequations}

With either representation for the stochastic gradient \eqref{e:nablaObjElig} or \eqref{e:nablaObjQ} we obtain an asymptotically unbiased SGD algorithm using $ \SGDnabla_\Obj(n) = \SGDnabla_\Obj^{\theta_n}$. Each of those described here are episodic: data is collected over the period $  0\le k\le \tstop(n)$ with $\theta_n$ fixed, and the input defined using $\feex^{\theta_n}$.

The \textit{natural gradient descent} algorithm updates the matrix gain via  $G_n = \haR_n^{-1}$,
where $\haR_0 >0$ with updates obtained recursively,
\begin{equation}
	\begin{aligned}
		\haR_{n}  &= \haR_{n-1}  + \stepf_{n}  [ - \haR_{n-1} +R_n  ] \,,\quad n\ge 1\,,
		\\
		R_n   &  = \sum_{k=0}^ {\tstop(n) } \logderfeex^{\theta_n}_k  [\logderfeex^{\theta_n}_k]^\transpose
		\,, \qquad   \logderfeex^{\theta_n}_k = \logderfeex^{\theta_n} ( \InfoState_k ,U_k)    
	\end{aligned} 
	\label{e:NGD}
\end{equation}
with $\stepf_n\gg \alpha_n$   (see \cite[Ch.~10]{CSRL}).

The two representations for the gradients prompt two choices for the stochastic gradient.   We focus here on the first,
\[
\SGDnabla_\Obj(n)  =  \sum_{k=0}^{\tstop(n) }  c(\Phi_k, U_k  ) S_k^{\theta_n}   
\]
leaving out the extension of the standard algorithm based on TD(1) learning to estimate the Q-function \cite[Ch.~10]{CSRL}.


The Polyak-Ruppert (PR) estimates are defined by
\begin{equation}
\thetaPR_n = \frac{1}{n} \sum_{i=1}^n  \theta_i\,, \qquad n\ge 1\, .
\label{e:PR}
\end{equation}
 Its \textit{asymptotic covariance} is defined as 
 \begin{equation}
\SigmaPR =
\lim_{n\to\infty} n \, \Expect[ \tilthetaPR_n \{\tilthetaPR_n  \}^\transpose]   
\label{e:PRcov}
\end{equation}
 When this exists and is finite, then the estimates 
 achieve the optimal mean-square convergence rate of $O(1/n)$.

%

The following is a consequence of recent stochastic approximation theory in \cite{chedevborkonmey21}.
 
 \begin{proposition}
\label[proposition]{t:SGDconverges}
Suppose that the assumptions of \Cref{t:ScoreRep} hold,  and in addition  
{(i)}
$\Obj$ is coercive with unique minimum $\theta^*$ and 
$\nabla\Obj$ is globally Lipschitz continuous.
{(ii)}
 $A^*\eqdef \nabla^2\Obj\, (\theta^*) $ is Hurwitz,   and the steady-state  covariance $R^* = \Cov( \logderfeex^{\theta^*} )$ is full rank.
(iii) The stepsize sequence is $\alpha_n = \alpha_0 n^{-\rho}$ with $1/2 <\rho <1$ and $\alpha_0>0$.

Then,   the SGD algorithm \eqref{e:QCDAC} is convergent almost surely and in mean square.  The PR-estimates are also convergent in both senses. 

The Central Limit Theorem (CLT)  holds, as well as the limit \eqref{e:PRcov},
 in which the asymptotic covariance is     $\SigmaPR =   [ ( R^* )^{-1}A^* ]^\transpose \Sigma_\nabla^* ( R^* )^{-1} A^*$  
 with
  $\Sigma_\nabla^*$ is the steady-state covariance of 
\eqref{e:nablaObjQ}  using the  policy $   \feex^{\theta^*}$.
 \qed
\end{proposition}

\wham{Example} 
Consider the one-dimensional family of policies in which $\theta$ approximates a threshold rule:  for a fixed large constant $\invTemp>0$, define 
$\feex^\theta(u\mid w)  =  [1+  \exp(\invTemp   [w-\theta] )]^{-1}  \exp(\invTemp u [w-\theta] )$,   so that the score function is 
\[
\logderfeex^\theta (u\mid w) = -\invTemp u   + \invTemp \feex^\theta(1\mid w)     
\]
In this scalar example we can adapt the natural gradient actor critic method to estimate $\nabla \Obj\, (\theta)$ for any fixed $\theta$. 

   \begin{figure}[h]
	\centering		 \includegraphics[width=1\hsize]{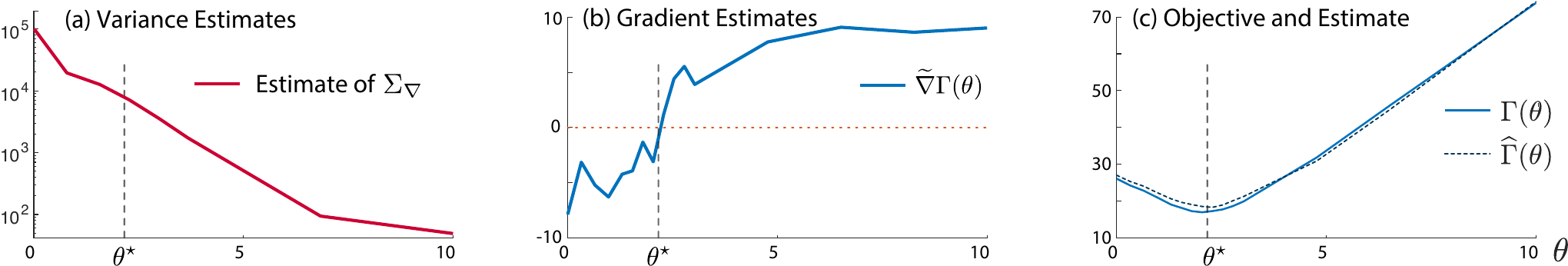}
	
	\caption{Statistics of the gradient estimate as a function of $\theta$,   based on the 
		Actor Critic method:
		(a)   Empirical variance of the gradient estimates.
		(b)  Gradient estimates using $N=10^4$ episodes.
		(c)   Objective and its approximation obtained from integrating the gradient estimate.				}
	\label[figure]{f:AC}
\end{figure}

\Cref{f:AC} shows results from a typical experiment using $\invTemp=20$.   Details on the simulation environment designed to obtain these approximations are postponed to  the Appendix.

 Rather than demonstrate results from an application of SGD,  the first two plots show estimates of the 
 mean 
 and variance  of the random variable in the expectation \eqref{e:nablaObjQ} for a range of values of $\theta$;  
 the precise means are  $\Sigma_\nabla \,(\theta) $ and  $  \nabla\Obj\, (\theta)$.

 \Cref{f:AC}~(b) shows gradient estimates, and  \Cref{f:AC}~(c) compares estimates of the objective function  $\Obj(\theta)$ obtained via standard Monte-Carlo, and the estimate obtained from integrating the  gradient estimates in (b):      $\widehat\Obj (\theta)\eqdef  \kappa + \int_0^\theta  \widehat\nabla\Obj (r)\, dr$. 
 The results indicates good news:  in spite of the enormous variance shown in (a),  especially large for smaller values of $\theta$,  
 the zero of the gradient estimate  $  \widehat\nabla\Obj\, (\theta)$  is very close to the optimal threshold value for CUSUM.   
 However, the massive variance presents a challenge in running the actor-critic algorithm to estimate $\theta^*$.

\wham{Q-learning}

Recall the solution to  the POMDP model in which the optimal policy is a function of an information state.
Consider  the canonical example in which this is the belief state (the sequence of conditional distributions $\{\condDist_k\}$) evolving on the unit simplex $\clS$,  and assume   that the underlying Markov chain $\bfPhi$ evolves on a finite set so that the simplex is finite-dimensional.

 The Q-function $Q^* \colon \clS \times\ustate\to\Re$ is the optimal value function associated with the objective \eqref{e:ObjPOMDPQCD}.    To place the equations in standard form denote $c(x,u) = (1-u) c_\circ(x) + u c_\bullet(x)$ for $x\in\state$ and $u\in \{0,1\}$ (recall   \eqref{e:ObjPOMDPQCD}).   
For any $u$ and $\upbeta\in\clS$  denote 
$  \clC(\upbeta,u) = \sum_x\upbeta(x)  c(x,u) $.    

 The value $Q^*(\upbeta,u)$ is defined to be  the minimum of 
$  \sum_\phi \upbeta(\phi) J(\phi , U_0^\infty ) $ over all admissible $U_1^\infty$,  subject to $(\condDist_0,U_0) = (\upbeta,u)$.  
It satisfies the dynamic programming (DP) equation,
\[
Q^*(\upbeta,u)  = \clC(\upbeta,u)    
+   \Expect[ \uQ^*( \condDist_{k+1} )  \mid  \condDist_k=\upbeta,U_k =u  ]  
\]
with   
$\uH(\upbeta) = \min \{ H(\upbeta, 0), H(\upbeta, 1) \}$ for any function $H\colon \clS \times \ustate\to\Re$.

Q-learning algorithms are based on the  characterization:     $ \Expect  \bigl[ \Tdiff^*_{k+1}   \mid \clY_k\big] =0$  or each $k$ and   any adapted input, with $  \Tdiff^*_{k+1}  =  - Q^*( \condDist_k , U_k)   +c_k
 							+ \uQ^*( \condDist_{k+1} )  $.   with $c_k = (1-U_k)\ind \{\tchange < k\} +  \kappa \,  U_k   (\tchange-k)_+$.

This motivates typical Q-learning algorithms.

Given a parameterized family of real-valued functions $\{Q^\theta : \theta\in\Re^d\}$ on $\sstate\times \ustate$,     the goal is to solve the \textit{projected Bellman equation}:  $\barf(\theta^*) =0$ with
\begin{equation}
\barf(\theta) \eqdef \Expect\bigl[ \{  
 - Q^\theta( \InfoState_k , U_k)   + c_k   
 							+ \uQ^\theta( \InfoState_{k+1} ) \}  \zeta_k  \bigr]
\label{e:projBE}
\end{equation}
where $\{\zeta_k \}$ is a $d$-dimensional stochastic process adapted to the observations.    

It is typical to take $\zeta_k = \nabla_\theta  Q^\theta( \InfoState_k , U_k) \big|_{\theta_k}$,  with $\theta_k$ the estimate at iteration $k$.    Theory to-date is largely restricted to a linear function class in which $Q^\theta =\theta^\transpose \psi$ with $\psi\colon \sstate\times \ustate  \to\Re^d$, and in this case $ \zeta_k = \psi( \InfoState_k , U_k) $.

Data required in an algorithm is based on successive runs up to time  $\tstop(n)$ for $n\ge 1$, and with $\tstop(0)\eqdef 0$,  which results in the observations $\{\InfoState_k^n ,\InfoState_{k+1}^n, U_k^n  ,c(\Phi_k^n, U_k^n)   \}$ for   $0\le k<\tstop(n)$.   
We suppress dependency on $n$ by stringing data together,  so that for example
\[
U_k  \eqdef  U_{ k- \tstop(n-1)}^n \  \ \textit{for $ \tstop(n-1) \le k < \tstop(n)$ and $n\ge 1$}
\]
with $\tstop(0)\eqdef 0$.

\begin{subequations}

A version of Q-learning is expressed as the recursion 
\begin{align}
  \theta_{k+1} & = \theta_k  + \alpha_{k+1}  G_k \zeta_k \Tdiff_{k+1}   \,,\qquad\qquad  k \ge 0
\label{e:Q-recursion}
   \\ 
   \Tdiff_{k+1} & =  - Q^{\theta_k} ( \InfoState_k , U_k)   + c_k 
 							+ \uQ^{\theta_k} ( \InfoState_{k+1} )   
\end{align} 
where the matrix gain sequence $\{G_k\}$ is a design choice;    
Zap Q-learning is in some sense \textit{optimal} \cite{CSRL}; 
This matrix gain was   used  in    \cite{chedevbusmey20a} for applications to optimal stopping.

\end{subequations}

We cannot apply   \cite{chedevbusmey20a}, based on the elegant algorithm of  \cite{tsivan99},  since the resulting policy will depend on the cost  $\{  c(\Phi_k, U_k)  \}$
(assumed observed in this prior work).

While there is great empirical success in the history of Q-learning,  to-date we only have  general conditions for stability of the algorithm, and existence of a solution to the projected Bellman equation \cite{mey23}.  
 Most crucial is the requirement that the input used for training is an $\epsy$-greedy policy (or a smoothed variant).
 It is shown that, subject to a mild full rank condition for $\psi$,  that for sufficiently small $\epsy>0$  the algorithm \eqref{e:Q-recursion}  is stable in the sense of ultimate boundedness,  and there exists at least one solution $\theta^*\in\Re^d$ to the projected Bellman equation  $\barf(\theta^*) =0$.   Convergence remains a topic of research.

Stability of Zap Q-learning with an oblivious policy (independent of parameter) is virtually universal
\cite{chedevbusmey20b}, but this paper makes no claims of existence of $\theta^*$ in this setting.    It is very likely that the main result of \cite{chedevbusmey20b} can be extended to $\epsy$-greedy policies.

\section{Numerical Results with Q-learning} 
\label{s:Qqcd}

This section contains two subsections:  the first illustrates the large-$\kappa$ approximations of \Cref{t:barcApprox}, and the second summarizes results obtained using  the Q-learning formulations described in \Cref{s:RLQCD}.

\wham{QCD model}

The conditional i.i.d.\ model  \eqref{e:QCDmodel} was used to generate observations, in which  $\preObs_k \sim N(0,\sigma^2) = \preDens$ and 
$\postObs_k \sim N(\mu_1,\sigma^2)  = \postDens$ 
with $\mu_1 = 0.5$ and $\sigma = 1$.

\begin{wrapfigure}[11]{r}{0.34\hsize}
	\vspace{-2.3em}					 
	\centering     
	\includegraphics[width=\hsize]{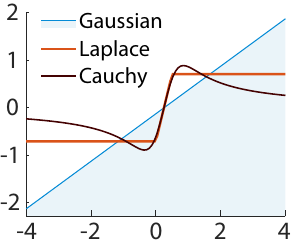}
	\vspace*{-5.5mm} 
	\caption{Three LLRs}
	\label[figure]{f:3LLRs}
\end{wrapfigure}
\smallskip

Three choices of   $F$ were tested, each of the form $F=\log(\surf_1/\surf_0)$:
\\
\textbf{Case 1:}   The ideal
Gaussian case, in which $ \surf_1 = \postDens$ and $\surf_0   =\preDens$. 
\\ \textbf{Case 2:}   $\surf_0$ is Laplace$(0,b)$ and $\surf_1 $ Laplace$(\mu_1,b)$ with $\mu_1 = 0.5$ and $b = \sqrt{\sigma^2/2}$ (matching second order statistics).
\\
\textbf{Case 3:} $\surf_0$ is Cauchy$(0,\disc)$ and $\surf_1 $ Cauchy$(x_1,\disc)$ with $x_1 = 0.5$ and $\CauchyScale$ chosen so that the Gaussian and Cauchy cdfs evaluated at $\sigma = 1$ are equal.

\begin{figure*}
	\centering		
	\includegraphics[width=0.995\hsize]{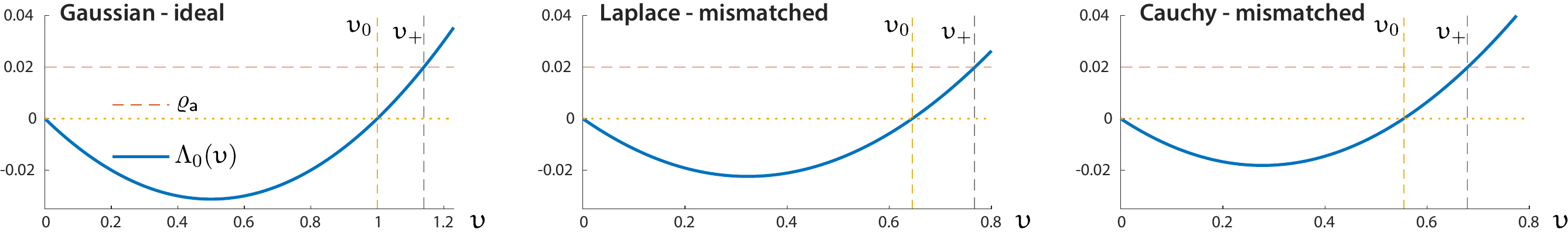}
	\caption{$\Lambda_0(\thexp)$ for ideal Gaussian alongside Laplace and Cauchy mismatched detectors}
	\label[figure]{f:Lam0Plots}
\end{figure*}

Plots of the three functions are shown in     \Cref{f:3LLRs},  and the corresponding log moment generating functions are shown in   \Cref{f:Lam0Plots}.

\subsection{Cost approximation}
\label{s:costapprox}

For the ideal case $F$ is the true LLR  for the marginals of  $\bfpreObs, \bfpostObs$, giving  $F(x)=L(x)=\mu_1 x - \mu_1^2/2$ and  $
\Lambda_0(\thexp) = m_1 \thexp  (\thexp -1)$. 
Hence the equation $\expa = \Lambda_0(\thexp_+)$ is easily inverted to obtain $\thexp_+ > 0$ (recall this is required in \eqref{e:ApproxJH}).   For the mismatched case,   $\thexp_+$ was approximated based on numerical computation of $\Lambda_0$.

\begin{subequations}
	For each choice of $F$  we obtained estimates via Monte-Carlo of  the optimal threshold and corresponding total cost defined in \eqref{e:optThresholdCost},  denoted $\{   \hathresh^*(\kappa),  \haJ^*(\kappa)  \}$.      In the plots described next the shifted values are displayed: 
	\begin{align}
		\barthresh_{\tiny{\sf s}}(\kappa)  &=
		\barthresh_\infty^*(\kappa) - \barthresh_\infty^*(100) + \hathresh^*(100)
		\label{s:shiftH}
		\\
		\barJ_{\tiny{\sf s}}(\kappa)  &=  \barJ_\infty^*(\kappa)  -
		\barJ_\infty^*(100)    +  \haJ^*(100)
		\label{s:shiftJ}
	\end{align}
	\label{s:shiftHJ}This ensures that the approximations \eqref{e:ApproxJH} coincide with the Monte-Carlo estimates of  \eqref{e:optThresholdCost} at $\kappa=100$.
\end{subequations}

\Cref{f:CostAndThresholdApproximations} shows the Monte Carlo estimates with $1\sigma$ confidence intervals compared to results from experiments in all three cases,   and two  different choices for the change time, each satisfying
$\expa=0.02$.  
In the first column the change time has geometric distribution with parameter $\expa$,  expressed $\tchange \sim \text{geo}(\expa)$.
In the second column the distribution is a mixture of geometrics:  $\tchange \sim \text{geo}(\expa)$ with probability $0.05$; else $\tchange \sim \text{geo}(0.2)$.      

The approximations are remarkably accurate in all cases,  which means that    the error $ \barthresh_\infty^*(\kappa) - \thresh^*(\kappa) $ is nearly constant over the entire range.  
Unfortunately, the constant value is large, which motivates learning techniques to obtain a near-optimal threshold.

\subsection{Q-learning}
The remainder of this section presents the design and evaluation of Q-learning for this Bayesian QCD problem.    
In each experiments the observations come from Shiryaev's model in which the change time is geometrically distributed,  and hence  the true optimal test is available.

\wham{Basis selection}
The basis for the function class $\{Q^\theta =\theta^\transpose \psi\}$  took the form $\psi(x,u) = (1-u) \psi^0(x)  + u \psi^1(x)$ (recall \Cref{s:RLQCD}).   
A four dimensional basis gave good results:  
\begin{align}
	&\psi^0(x)  = [x; q(x); 0;0;0] \nonumber
	\\
	&\psi^1(x) =  [0;0; 1;x; q(x)) ] 
	\label{e:basis2}
\end{align} 
with $q(x) =  x\exp(-x/b_q)$ for a choice of constant $b_q$.

\begin{figure}[h]
	\centering		
	\includegraphics[width=0.8\hsize]{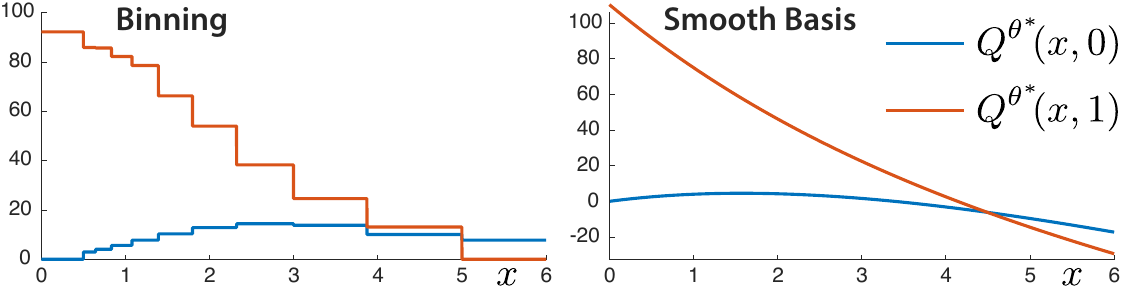}
	\caption{Insights from binning led to the basis in \eqref{e:basis2}}
	\label[figure]{f:Qsolutions}
\end{figure}

This basis was chosen based on preliminary experiments with a particular choice of binning:
$\psi_i(x,u) = \ind\{x\in S_{k_i} ,   u = u^i\} $  for a collection of intervals $ \{S_j \}$ and input values $\{u^j\} $.   Further details are provided below;  see in particular \Cref{f:Qsolutions}
and surrounding discussion.

\wham{Stability of Q-learning}

\begin{figure}[h]
	\centering		\includegraphics[width=.9\hsize]{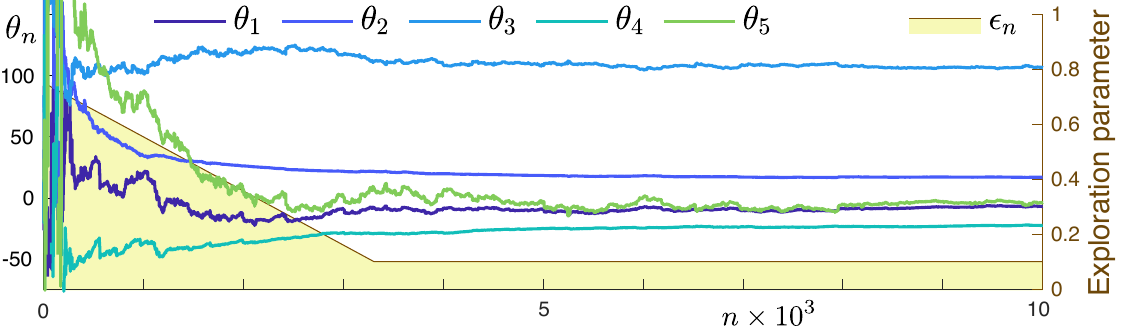}
	\caption{Parameter estimates using a decaying exploration schedule.}  
	\label[figure]{f:TypicalTrajectories}
\end{figure}

Recent theory recalled in \Cref{s:RLQCD} shows that exploration implies stability of Q-learning under mild assumptions on the basis and the  oblivious policy.
Crucial for stability is that the exploration gain be sufficiently small.   In the experiments surveyed here this was taken  to be time varying , with a typical choice illustrated in  \Cref{f:TypicalTrajectories}:
$\epsy_n =  \max\{\epsy_f ,  \epsy_0 +  (1 - n/n_0) (\epsy_f - \epsy_0)\}$,
defined so that $\epsy_n = \epsy_f < \epsy_0 $ for $n\ge n_0$. 
Consistent with theory from \cite{mey23}, it was found in experiments that both  the scalar gain algorithm and Zap Q-learning were convergent for sufficiently small $\epsy_f$.

In the case of Zap Q-learning, convergence held for a wide range of $\epsy_f$;   the value $\epsy_f =0.1$ was used in the experiments described in the following.   The scalar gain algorithm was \textit{far less reliable}:  
$\epsy_f\le 10^{-4}$ was required for stability.  

One can then ask, is the  parameter estimate $\theta^*$ obtained using Zap Q-learning consistent with the output of the scalar gain algorithm?   Let $\barf$ denote the vector field for the mean flow associated with the scalar gain Q-learning algorithm---see \eqref{e:projBE}.     We computed the Jacobian  $\barA(\theta^*) = \partial_\theta \barf (\theta^*)$ for the parameter $\theta^*$ obtained using Zap Q-learning, and discovered that in most cases it had at least one eigenvalue in the strict right half plane in $\Co$.     In such cases, $\theta^*$ is not asymptotically stable for the mean flow.   Standard stochastic approximation theory implies that the parameter  $\theta^*$ would not be found using the Q-learning  algorithm \eqref{e:Q-recursion} with $G_k=I$. 



\wham{Oblivious policy}
Exploration was designed to depend on $\kappa$:  at the start of episode $i$,   a threshold $\thresh^{\epsy,i}(\kappa)$  was drawn uniformly at random from an interval $[a_{\kappa},b_{\kappa}]$.    Then, $U_n = \ind \{ \InfoState_n \ge \thresh^{\kappa,i} \}$ for each $n$ in this episode. 

The threshold approximations of \Cref{e:ApproxJH}  motivated the design of the interval:  
\begin{equation}
	[a_{\kappa},b_{\kappa}] =     [   \barthresh_\infty^*( \kappa) + \eta   -   \delta,    \barthresh_\infty^*( \kappa) + \eta   +   \delta ]
	\label{e:abinterval}
\end{equation}
with $\delta>\eta>0$ constant.
This ensured significant exploration in all cases considered. 


\wham{Numerical experiments}

Each parameter $\theta\in\Re^d$ defines a policy,   $\fee^\theta(x) = \argmin_u Q^\theta(x, u)$ for $x\in\Re_+$.
In the applications considered here this becomes
\begin{equation}
	\fee^\theta(x) =  \ind\{ Q^\theta(x, 0) \ge Q^\theta(x, 1) \}     \,,  \qquad x\in\Re_+\, .
	\label{e:fee_theta_HTorQCD}
\end{equation}
In every successful application of Q-learning it was found that this policy had a threshold form
\begin{equation}
	\fee^\theta(x) =   \ind\{x\ge \thresh^\theta\}\,, \quad \thresh^\theta >0
	\label{e:fee_theta_threshold}
\end{equation}

\begin{wrapfigure}[15]{r}{0.35\hsize}
	\vspace{-0.5em}	 
	\centering   
	\includegraphics[width=\hsize]{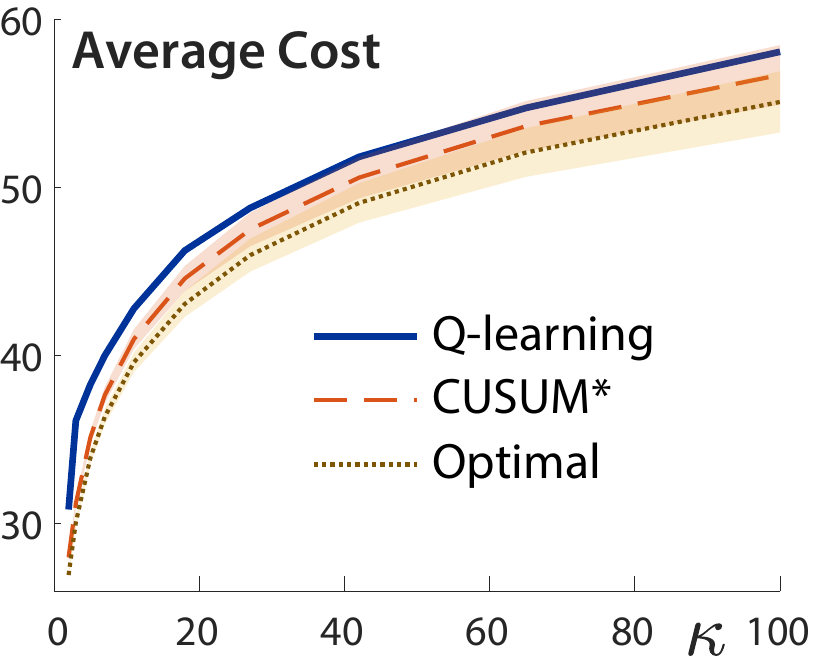}
	\vspace*{-5.5mm} 
	\caption{Average cost comparisons.}    
	\label[figure]{f:Performance}
\end{wrapfigure}

Algorithm performance is investigated in the remainder of this section.    For each algorithm, PR-averaging was used to define the final estimate $\hatheta$,   and from this a final policy $ \hafee \eqdef \fee^{\hatheta}$ whose performance is compared to the optimal.


Initial experiments involved a choice of binning,  resulting in $d=2(d_0-1)$, with $d_0$ the number of bins. 
However, binning proved insufficient for obtaining  thresholds close to optimal over all $\kappa$, due in part to a dependence on the choice of bin spacing. This shortcoming is illustrated in \Cref{f:Qsolutions}, where bin spacing influences the intersection $Q^\theta(x,0)=Q^\theta(x,1)$ for the policy  \eqref{e:fee_theta_HTorQCD}. This inspired the ``smooth" basis in \eqref{e:basis2} for which we observed two advantages compared to binning: 1/ better performance of $\hafee$ for all $\kappa$, and 2/ reduced training time.

Histograms were generated to evaluate the variance of the parameter estimates. Let $\EpsLength=\EpsLength(N)\ge N$ denote the total number of samples $(\InfoState_k, U_k)$    collected over $N$ episodes, so that  $\hatheta = \thetaPR_\EpsLength$ is the final estimate.

The \textit{asymptotic covariance} of   $\tilthetaPR_\EpsLength \eqdef \thetaPR_\EpsLength-\theta^*$ is denoted
\[
\SigmaPR =
\lim_{N\to\infty} \EpsLength \Expect[  \, \tilthetaPR_\EpsLength \{\tilthetaPR_\EpsLength  \}^\transpose]
\]
This was estimated using the \textit{batch means method}:    $M$ independent runs resulted in the estimates  $\{ \thetaPR_{\EpsLength^i} ,  \EpsLength^i: 1\le i\le M\}$.  The empirical covariance of $\{ Z^i = \sqrt{\EpsLength^i} [ \thetaPR_{\EpsLength^i} - \barthetaPR ] : 1\le i\le M \}$ provides an estimate of
$\SigmaPR$ with vanishing error as $M\to\infty$ and then $ N\to\infty$.

\begin{figure}[h]
	\centering		 
	\includegraphics[width=0.75\hsize]{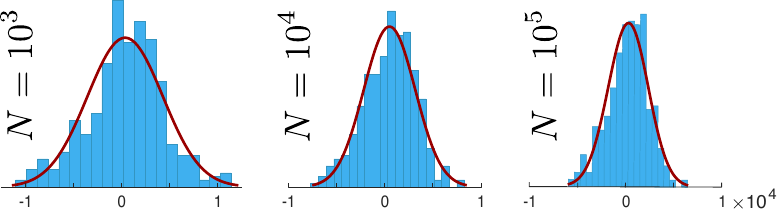}
	\caption{Histograms of $\{Z^i_1 : 1\le i\le M\}$  for three values of $N$. }
	\label[figure]{f:Hists}
\end{figure}

An example is shown in \Cref{f:Hists} for the case $\kappa=27$,  using $M=400$ and three different values of $N$.   
Only the fourth component of the five dimensional histogram is shown, giving an estimate of $\SigmaPR(4,4)$.
What is crucial here is that the estimate of this value 
is nearly identical for the three values of $N$ chosen.     Similar results were observed for estimates the other diagonal entries of $\SigmaPR$.

%

This is an example of how the CLT can be used to estimate required run lengths by first conducting a large number of independent experiments with a relatively short run length---in this example,  
$N=10^4$ provides a reasonable estimate of the variance of $Z^i = \sqrt{\EpsLength^i} [ \thetaPR_{\EpsLength^i} - \barthetaPR ]$ for each $i$ and $N\gg 10^4$.

\begin{figure}[h]
	\centering	 
	\includegraphics[width=0.75\hsize]{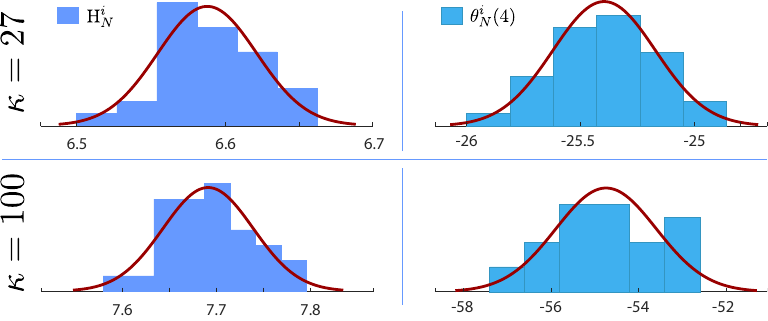}
	\caption{Histograms of $\{ \thresh_N^{\theta,i} ,  \theta_N^i(4): 1\le i\le M\}$ for $N=10^6$}		
	\label[figure]{f:HistsPar+Thresh27+100}
\end{figure}

Thresholds $\thresh_N^{\theta,i}$ yielded smaller empirical variance than $\theta_N^i$ for all $\kappa$. Example histograms showing the fourth parameter $\theta_N^i(4)$ are included in \Cref{f:HistsPar+Thresh27+100} for $M=30$ and $N=10^6$.

\Cref{f:Performance} shows the average cost of the policy $\hafee$ for the ideal Gaussian case. Also included are estimates $\haJ^*(\kappa)$ using CUSUM* and the \textit{optimal test}---Shiryaev test defined below  \eqref{e:QCDmodel}.   The 
 confidence intervals for each are of $1\sigma$ standard deviation. 
One metric for success is that the shape of the average cost curve obtained through Q-learning resembles the optimal, which we observed for $\kappa \ge 20$.

\whamit{Mismatched cases.} 

Experiments were repeated for Cases 2 and 3: the surrogate information state  $\{\InfoState_n\} $ differs in each case based on the respective LLRs, plotted in \Cref{f:3LLRs}.  

Recall the discussion surrounding \eqref{e:fee_theta_HTorQCD}:  the policy $\fee^\theta$ obtained from Q-learning resulted in a threshold policy in almost all cases.    Therefore, the best performance from Q-learning can be no better than CUSUM*,  which uses the optimal threshold 	\eqref{e:optThreshold}.   Moreover, the performance in Cases 2 and 3 is poor when compared with   Case 1---recall the comparisons in  \Cref{f:CostAndThresholdApproximations}.

However, in all three cases, the policy obtained from Q-learning using $\kappa \ge 27$ yielded  average cost within 5\% of its respective CUSUM* cost estimate $\haJ^*(\kappa)$.

\section{Conclusions} 
\label{s:conc}

The theory and numerical results in this paper motivate many directions for future research:

\whamb   The preliminary  research surveyed in \Cref{s:RLQCD} motivates research 
on variance reduction for  stochastic gradient descent.

\whamb     The independence assumptions in \Cref{t:barcApprox} are unfortunate.     Relaxation of the i.i.d.\ assumption is the main topic of \cite{coomey24b,coomey24c},   along with conditions under which the  independence of $\tchange$ and $\bfmX^0$ may be relaxed;
the assumptions of \Cref{t:POMDPhazard} are sufficient  to obtain an extension of \Cref{t:barcApprox}.

\whamb
In applications of interest to us there may be well understood behavior before a change  (which might represent a fault in a transmission line,  or a computer attack).    We cannot expect to have  a full understanding of post-change behavior.  The choice of surrogate information state must be reconsidered in these settings,  perhaps based on techniques for universal hypothesis testing     (see \cite{huaunnmeyveesur11,huamey13} and the references therein).

%
%
%

 \bibliographystyle{abbrv}
 \bibliography{arXiv_QCD_Irma}

 \appendix

\section{Appendix}

 This Appendix contains details on the QCD numerical results.

\wham{Actor-critic method}

Estimates of $\Sigma_\nabla$ were obtained by averaging $N=10^7$ independent episodes. 

Estimates of 
the gradient   $\SGDnabla_\Obj(\theta)$ were obtained using a much shorter run,  with $N=10^4$.
Two estimates of the objective were produced with $N=10^4$:  $ \Obj(\theta) $ using Monte-Carlo to estimate the expectation in \eqref{e:Obj}
directly and $\widehat\Obj (\theta)$ through numerical integration.

\wham{Q-learning for QCD} 

For the ideal and mismatched cases, Monte Carlo simulations were used to estimate $\MDE$ and $\MDD$ for CUSUM*. Estimates were also obtained for the optimal Shiryaev test $\InfoState_n =p_n=  \Prob\{  \tchange \le  n \mid \Obs_0^n \}.$ Parameters for the stochastic processes $\{\preObs_k , \postObs_k$\} and $\tchange$ matched those used for Q-learning. For both simulations, $N=2e4$ sample paths were run. A range of $T = 10^3$ thresholds $0 \le \thresh \le 20$ was used for CUSUM*, and $0 \le \thresh \le 1$ for Shiryaev. For each $\thresh$, a pair $\MDE(\thresh)$ and $\MDD(\thresh)$ was obtained by averaging the eagerness and delay over $N$ runs. This repeated for $M=200$ independent runs, averaging again to obtain for each test a $T \times 2$ matrix $[\MDE, \MDD]$, where each row corresponds to a different threshold $\thresh$. Estimates of these quantities are random variables, whose variances were found to be very small. Estimates $ \hathresh^*(\kappa)$ and  $\haJ^*(\kappa) $ as described before \Cref{s:shiftHJ} were then obtained for a range $2 \le \kappa \le 100$ to generate the average curves in \Cref{f:Performance}. 

Additional elements of the Q-learning experiments include:

\wham{1.)}   \textit{Parameter initializations.} 
In Q-learning,  for each independent episode, $\theta_0$ was chosen uniformly at random in [-100,100].

\wham{2.)}  \textit{Resetting.} 
It is typical in any stochastic approximation implementation to observe large transients.   In these experiments,    $\theta_{k+1}$ was sampled independently of the past, 
uniformly at random in [-100,100],  if  $\Vert \theta_k \Vert_\infty > 5e3$.



\end{document}